\documentclass[10pt]{amsart}
\usepackage{amsmath,amscd}
\usepackage{amssymb,xypic,graphicx}
\usepackage{amsthm,verbatim}

\usepackage{hyperref}
\usepackage{tikz}
\usepackage{tikz-cd}
\newcommand*{\rom}[1]{\expandafter\@slowromancap\romannumeral #1@}
\newcommand{\tr}{{\sf tr}}

\newcommand{\fkg}{{\mathfrak{{g}}}}

\newcommand{\irr}{{\sf irr}}

\newcommand{\virt}{{\sf virt}}

\newcommand{\Cone}{\operatorname{Cone}}

\newcommand{\Sym}{\operatorname{Sym}}

\newcommand{\KK}{{\mathcal{K}}}

\newcommand{\II}{{\mathcal I}}

\newcommand{\BB}{{\mathcal{B}}}

\numberwithin{equation}{subsection}

\newtheorem{thm}{Theorem}[section]
\newtheorem{prop}[thm]{Proposition}
\newtheorem{lem}[thm]{Lemma}

{  \theoremstyle{definition}
\newtheorem{defi}[thm]{Definition}

\newtheorem{rem}[thm]{Remark}

}

\newcommand{\id}{\operatorname{id}}

\newcommand{\ra}{\rightarrow}
\renewcommand{\AA}{{\mathcal A}}

\newcommand{\TT}{{\mathcal T}}

\newcommand{\Om}{\Omega}

\newcommand{\Hom}{\operatorname{Hom}}

\newcommand{\R}{{\Bbb R}}
\newcommand{\Z}{{\Bbb Z}}

\newcommand{\ot}{\otimes}

\newcommand{\ed}{\qed\vspace{3mm}}

\newcommand{\SU}{{\sf SU}}
\newcommand{\su}{{\sf su}}

\newcommand{\MC}{\sf MC}

\newcommand{\cob}{\sf cob}

\numberwithin{equation}{section}

\title{Casson invariants via virtual counting}
\author{Junwu Tu}
\thanks{Mathematics Department, University of Missouri, Columbia, MO, 65211, USA. E-mail: {\texttt{ tuju@missouri.edu}}}

\begin{document}
\begin{abstract}
This is a research announcement on an alternative definition of the Casson invariants by means of virtual counting of the moduli space of irreducible representations of the fundamental group into $\SU(2)$. Along the way, by using derived differential geometry, we propose a general framework to obtain invariants from Chern-Simons type gauge theories. 
\end{abstract}

\maketitle

\section{Introduction} The Casson invariant is a differential invariant associated to compact oriented integral homology $3$-spheres. Let $M$ be such a manifold.  Denote by 
\[ R^\irr(M):=\Hom^\irr(\pi_1(M),\SU(2))/\SU(2),\]
the moduli space of irreducible representations of its fundamental group into $\SU(2)$. Recall that a representation is called irreducible if its stabilizer group is equal to $\Z/2\Z$, the center of $\SU(2)$. The integral sphere condition implies that the class of the trivial representation is the only reducible representation. Furthermore, it is an isolated point in $R(M)$, the moduli space of all representations. Since $\SU(2)$ is compact and $\pi_1(M)$ is finitely generated, the moduli space $R(M)$ is compact, which implies that the topological space $R^\irr(M)$ is also compact. The definition of Casson invaraints requires the choice of a Heegaard splitting $$M=H_1\cup_\Sigma H_2$$ with $H_1$, $H_2$ handlebodies of genus $g$ with boundary a closed surface $\Sigma$. Then the Seifert-Van Kampen theorem implies a fiber diagram of topological spaces
\begin{equation}\begin{CD}~\label{diagram:classical}
R^\irr(M) @>>> R^\irr(H_1)\\
@VVV               @VVV\\
R^\irr(H_2) @>>> R^\irr(\Sigma).
\end{CD}\end{equation}
Even though $R^\irr(M)$ is not a manifold, both $R^\irr(H_1)$ and $R^\irr(H_2)$ are submanifolds in an ambient smooth manifold $R^\irr(\Sigma)$. Their dimensions are
\[ \dim R^\irr(\Sigma)=6g-6; \;\; \dim R^\irr(H_1)=3g-3; \;\; \dim R^\irr(H_2)=3g-3.\]
Furthermore, these three manifolds have orientations induced from that of $M$. The compactness of $R^\irr(M)$ implies that a small perturbation of $R^\irr(H_1)$ will intersect $R^\irr(H_2)$ transversely along finitely many points. The Casson invariant is then defined to be
\[ \lambda(M):=\frac{1}{2} \cdot \big( \mbox{ intersection number of $R^\irr(H_1)$ and $R^\irr(H_2)$}\big).\]

Intuitively speaking, Casson invariants should be thought of as simply counting invariants of the moduli space $R^\irr(M)$. In fact, its holomorphic version the Donaldson-Thomas invariants were defined by virtual counting of moduli spaces of stable coherent sheaves~\cite{Tho}, using the virtual cycle construction in the algebraic setting by~\cite{BF} and~\cite{LT}. In this paper, using derived differential geometry, we give a new definition of the Casson invariant by means of virtual counting of the moduli space $R^\irr(M)$. Our starting point is the Chern-Simons gauge theoretic description of the moduli space $R^\irr(M)$, which is already used by Taubes~\cite{Tau} to give a definition of the Casson invariant by Morse theory on Banach manifolds. This point of view allows one to construct certain derived enhancement the moduli space $R^\irr(M)$. Then the main observation is that the derived moduli space automatically carries yet another structure, Kuranishi structure from which one can obtain the virtual fundamental class. 

\[ \mbox{Chern-Simons theory} \Rightarrow \mbox{$[0,1]$-type homotopy $L_\infty$ spaces} \Rightarrow \mbox{Kuranishi manifolds}.\]


\medskip
\noindent Global Kuranishi structures was pioneered by Fukaya-Oh-Ohta-Ono~\cite{FO},~\cite{FOOO}. In order to set a solid foundation for studying moduli spaces of pseudo-holomorphic curves, there have been an extensive amount of ongoing research in this direction~\cite{Joyce},~\cite{MW},~\cite{MW2},~\cite{Pardon},~\cite{Yang}. In this paper, we shall present yet another definition of Kuranishi manifolds. This new definition is quite similar to Joyce's definition~\cite{Joyce}. For this reason, I consider it as no more than a simplified version of Joyce's theory where we work with local charts in Euclidean spaces and trivial vector bundles. However, it is interesting how we arrived at this definition: by writing down an existing definition in the theory of $L_\infty$ spaces. This enables the usage of global Kuranishi theory to study Chern-Simons type quantum field theories since it is known how to produce (homotopy) $L_\infty$ spaces from them~\cite{Tu}. In the following sections, we describe in more detail how this proposal works in the case of classical Chern-Simons theory, and how Casson invariants appear as a result of integration.

\medskip
\noindent {\bf Acknowledgment.} I am grateful to Professor Joyce for pointing out an important mistake in defining Kuranishi manifolds. I also thank Lino Amorim for useful discussions regarding Joyce's Kuranishi theory.

\section{$[0,1]$-type homotopy L-infinity spaces}
Again, we let $M$ be a compact oriented integral homology $3$-sphere. By the Chern-Simons theory on $M$, we mean the differential graded Lie algebra
\[\mathfrak{g}:=\Omega^*(M)\otimes \su(2),\]
together with the symmetric invariant pairing
\[ \langle \alpha,\beta \rangle:=\int_M \tr(\alpha\wedge\beta).\] 
We recall some basic definitions. An element $b\in \mathfrak{g}^1$ is called a Maurer-Cartan element if it satisfies the Maurer-Cartan equation
\[ db+\frac{1}{2}[b,b]=0.\]
Given a Maurer-Cartan element $b$, we obtain a perturbed differential $d_b:=d+[b,-]$ on $\mathfrak{g}$. Denote its cohomology by $H_b^*(\mathfrak{g})$. Two Maurer-Cartan elements $b_1$ and $b_2$ are called gauge equivalent if there exists a smooth map $g:M\ra \SU(2)$ such that
\[ b_2= gb_1g^{-1}-dg\cdot g^{-1}.\]
Equivalent Maurer-Cartan elements have isomorphic cohomology $H_{b_1}^*(\mathfrak{g})\cong H_{b_2}^*(\mathfrak{g})$, via the adjoint action of $g$. The Maurer-Cartan moduli space $\MC(\mathfrak{g})$ is the set of equivalence classes of Maurer-Cartan elements.  Given a Maurer-Cartan element $b$, taking the holonomy representation of the flat connection $d+[b,-]$ sets up a well-known bijection between $\MC(\mathfrak{g})\cong R(M)$.

Let $[b]\in \MC(\fkg)$ be a Maurer-Cartan class corresponding to an irreducible representation in $R(M)$. This implies in particular there is no continuous family of automorphisms, which implies that
\[ H^0_b(\fkg)=0.\]
By Poincar\'{e} duality, the group $H^3_b(\fkg)=0$ as well. Thus we have
\[ H^*_b(\fkg)=H^1_b(\fkg)\oplus H^2_b(\fkg),\]
with the two groups dual to each other by Poincar\'{e} duality. The two cohomology groups can be used to give a deformation-obstruction description of the moduli space $R^\irr(M)$, detailed in the following proposition. Its proof follows from~\cite[Section 2]{Tu}.

\medskip
\begin{prop}~\label{prop}
There exists a $L_\infty$ algebras structure on $H^*_b(\fkg)$, with the only non-trivial structure maps given by
\[ l^b_k: \Sym^k H^1_b(\fkg) \ra H^2_b(\fkg), \;\; k\geq 2.\]
Furthermore, the series
\[ \sum_{k\geq 2} \frac{1}{k!} l^b_k(x,\cdots, x)\]
converges in an open neighborhood $V$ of the origin in the vector space $H^1_b(\fkg)$, which defines an analytic section $s: V \ra \underline{H^2_b(\fkg)}$ of the trivial bundle generated by $H^2_b(\fkg)$ over $V$. By possibly shrinking the open subset $V$, there exists a homeomorphism
\[ \psi: U\ra s^{-1}(0),\]
where $U$ is an open neighborhood of $[b]$ in $R^\irr(M)$.
\end{prop}

\begin{rem}
Using Poincar\'{e} duality, the bundle $\underline{H^2_b(\fkg)}$ can be identified as $\Omega_V$, the bundle of $1$-forms on $V$. Under this identification, one can show that the section $s$ is in fact the differential of a real analytic function $V$ defined by
\[ f(x):=\sum_{k\geq 2} \frac{1}{(k+1)!}\langle l^b_k(x,\cdots, x),x\rangle.\]
\end{rem}

The proposition above leads us to consider a certain enhancement of the topological space $R^\irr(M)$, encoding the $L_\infty$ algebras governing the local deformation theory. The local model of this enhancement is described in the following

\medskip
\begin{defi}~\label{defi-lie}
A $[0,1]$-type $L_\infty$ space is a triple $(V,\fkg,D)$, where
\begin{itemize}
\item[(a.)] $\fkg$ is a curved $L_\infty$ algebra bundle of the form
\[ \fkg=T_V[-1]\oplus E[-2].\]
Denote by $(C^*\fkg,Q)$ its Chevalley-Eilenberg algebra. By definition, as a graded algebra bundle, $C^*\fkg$ is given by $\hat{S}\Om_V\ot \Lambda E^\vee$ with $\Om_V$ in degree $0$ and $E^\vee$ in degree $-1$.
\item[(b.)] $D: C^*\fkg \ra \Om_V \ot C^*\fkg$ is a flat connection on the bundle $C^*\fkg$ of cohomological degree zero. Furthermore, we require that
\begin{itemize}
\item[--] it is a derivation, i.e. $D(ab)=Da\cdot b + (-1)^{|a|} a\cdot Db$.
\item[--] the component $\Om_V \ra \Om_V\ot 1$ of $D$ is given by $D(\alpha)=\alpha\ot 1$.
\item[--] the connection $D$ is compatible with the $L_\infty$ structure in the sense that $ [D,Q]=0$.
\end{itemize}
\end{itemize}
We define the virtual dimension of a $L_\infty$ space $(V,\fkg,D)$ by $\dim V - {\sf rank\;} E$. Let $U$ be a topological space. A $L_\infty$ enhancement of $U$ is a quadruple $(V,\fkg,D,\psi)$ such that $(V,\fkg,D)$ is a $L_\infty$ space, and $\psi: U\ra V$ is a homeomorphism onto the subspace $s^{-1}(0)$ where $s$ is the curvature term of $\fkg$ (which by definition is a section of the bundle $E$). We refer to $U$ as the footprint of the enhancement $(V,\fkg,D,\psi)$.
\end{defi}

\begin{rem}
The notion of ``$L_\infty$ space" was first suggested by Costello~\cite{Cos}. The definition above is slightly different from Costello's original one by allowing a curvature term to appear on the $L_\infty$ algebra bundle $\fkg$. 
\end{rem}

A homomorphism between two $L_\infty$ enhancements (of the same topological space $U$) $(V_\alpha, \fkg_\alpha, D_\alpha, \psi_\alpha)$ and $(V_\beta,\fkg_\beta,D_\beta,\psi_\beta)$ is a pair $(f,f^\sharp)$ where $f: V_\alpha\ra V_\beta$ is a smooth map of manifolds which commutes with the footprint maps $\psi_\alpha$ and $\psi_\beta$, and $f^\sharp: \fkg_\alpha \ra f^* \fkg_\beta$ is a $L_\infty$ homomorphism which intertwines with the connections $D_\alpha$ and $f^*D_\beta$ on the corresponding Chevalley-Eilenberg algebra bundles.

One can also define homotopies between homomorphisms, and higher homotopies. A homotopy $L_\infty$ enhancement of a topological space is, roughly speaking, given by local $L_\infty$ enhancements on a covering, together with coherent homotopy gluing data on multiple intersections. The constructions in~\cite{Tu} can be used to prove the following
\medskip
\begin{thm}~\label{thm-lie}
The moduli space $R^\irr(M)$ admits a natural $[0,1]$-type homotopy $L_\infty$ enhancement denoted by $\mathfrak{R}^\irr(M)$.  Its virtual dimension is $0$.
\end{thm}

\begin{rem}
The results in this section do not need the assumption that $M$ is an integral homology $3$-sphere, nor the particular choice of the group $\SU(2)$. However, without this assumption the underlying moduli space $R^\irr(M)$ may not be compact. This is the key difficulty in defining Casson type invariants in the general setup.
\end{rem}

\section{Kuranishi manifolds}

Let $U$ be a topological space. A {\sl Kuranishi chart} of $U$ is given by a quadruple $(V,E,s,\psi)$ where $V$ is a smooth manifold, $E$ a finite rank vector bundle on $V$, $s$ a section of $E$ over $V$, and $$\psi: U\ra V$$ is a continuous map which gives a homeomorphism onto the subspace $s^{-1}(0)$. Again we call $U$ the footprint of the chart $(V,E,s,\psi)$.

Let $(V_\alpha,E_\alpha, s_\alpha, \psi_\alpha)$ and  $(V_\beta,E_\beta, s_\beta,\psi_\beta)$ be two Kuranishi charts with the same footprint $U$. A Kuranishi morphism 
is given by a pair $(f,\hat{f})$ where $f: V_\alpha\ra V_\beta$ is a smooth map such that $f\circ \psi_\alpha=\psi_\beta$, and $\hat{f}: E_\alpha \ra f^*E_\beta$ a bundle map such that $\hat{f}\circ s_\alpha= f^*s_\beta$.

The following extremely simple observation is the starting point of our take to build a global theory of Kuranishi spaces. The proof  is elementary which only requires unwinding the definitions.

\medskip
\begin{lem}~\label{intro-space}
Let $V\subset \R^n$ be an open subset. Let $\fkg:= T_V[-1]\oplus \underline{\R^m}[-2]$ be the graded bundle given the tangent bundle in degree $1$ and the trivial rank $m$ bundle in degree $2$. Fix a topological space $U$. Then there is an equivalence of categories:
\[  L_\infty \mbox{ enhancements $(V,\fkg,D,\psi)$ of } U \cong \mbox{Kuranishi charts $(V,\underline{\R^m},s,\psi)$ of } U\]
\end{lem}

The difficulty in obtaining a global Kuranishi theory lies at the delicate question of how to define coordinate changes between Kuranishi charts. In literature, there are quite a few similar, yet different definitions of this notion~\cite{FOOO},\cite{Joyce},\cite{MW},\cite{Pardon}. On the other hand,  in~\cite{Tu}, the notion of homotopy $L_\infty$ spaces which is a global version of $L_\infty$ spaces is already defined. In view of Lemma~\ref{intro-space}, a natural question is what if we translate the definition of homotopy $L_\infty$ spaces to the language in Kuranishi theory. After explicitly writing done this definition in terms of the tangent and obstruction bundles, we arrived at the following

\medskip
\begin{lem}~\label{intro-homotopy}
Let $(f_0,f_0^\sharp)$ and $(f_1,f_1^\sharp)$ be two morphisms between $L_\infty$ enhancements, and let $(f_0,\hat{f}_0)$, $(f_1,\hat{f}_1)$ be the associated morphism between Kuranishi charts. Then a homotopy in the sense of~\cite[Section 4.4]{Tu} is equivalent to the following data:
\begin{itemize}
\item[(a.)] A family of map
\[ F(t): V_\alpha \ra V_\beta,\; t\in [0,1]\]
such that $F(0)=f_0$, $F(1)=f_1$.
\item[(b.)] A family of bundle map $\widehat{F}(t): \underline{\R^{m_\alpha}} \ra F(t)^*\underline{\R^{m_\beta}}$ such that $\hat{F}(0)=\hat{f}_0$, $\widehat{F}(t)=\hat{f}_1$, and satisfy the equation
\[ F(t)^*s_\beta=\hat{F}(t) s_\alpha, \; \forall t\in [0,1].\]
\item[(c.)] Two families of bundle maps
\begin{align*}
\Lambda(t):& \underline{\R^{m_\alpha}} \ra F(t)^* \underline{\R^{n_\beta}},\\
\Xi(t):& \wedge^2\underline{\R^{m_\alpha}} \ra F(t)^*\underline{\R^{m_\beta}}.
\end{align*}
such that
\begin{align*}
\frac{dF(t)}{dt}|_x &= \Lambda(t)(s_\alpha(x)), \;\; \forall x\in V_\alpha\\
\frac{d\hat{F}(t)}{dt}|_{(x,\xi)} &= ds_\beta|_{F(t)(x)}\big(\Lambda(t)(x,\xi)\big)+\Xi(t)(s_\alpha(x)\wedge \xi),\;\; \forall (x,\xi)\in \underline{\R^{m_\alpha}}.
\end{align*}
Here in the second equation, $\frac{d\hat{F}(t)}{dt}|_{(x,\xi)}$ stands for the derivative valued in the fiber direction, and $ds_\beta$ is computed using the trivialization by considering $s_\beta$ as a map $V_\beta \ra \R^{m_\beta}$.
\end{itemize}
\end{lem}

The author tried to use Conditions $(a.)$,$(b.)$,$(c.)$ to define Kuranishi manifolds, but it turned out to be unnecessarily complicated as the required structures are often hard to construct. However, by working out some consequences of these conditions, we are led to the following simplified version with which we found most convenient to develop a global Kuranishi theory.

\medskip
\begin{defi}~\label{def:homotopy}
Two morphisms 
\[(f_0,\hat{f}_0), (f_1,\hat{f}_1): (V_\alpha,\underline{\R^{m_\alpha}} ,s_\alpha, \psi_\alpha) \ra (V_\beta,\underline{\R^{m_\beta}},s_\beta,\psi_\beta)\]
between Kuranishi charts of a topological space $U$ are called homotopic if there exists a quotient bundle map (see the footnote)
\[ \Lambda\in \KK\Hom(\underline{\R^{m_\alpha}} , \underline{\R^{n_\beta}})~\footnote{There are two subspaces of the space of bundle maps $\Hom(\underline{\R^{m_\alpha}}, \underline{\R^{n_\beta}})$ over $V_\alpha$:
\begin{align*}
K_1&:=\left\{ \Lambda\in \Hom(\underline{\R^{m_\alpha}}, \underline{\R^{n_\beta}})|\;\;\Lambda(s_\alpha)=0.\right\}\\
K_2&:=\left\{ \Lambda\in \Hom(\underline{\R^{m_\alpha}}, \underline{\R^{n_\beta}})|\;\; \Lambda|_U=0.\right\}
\end{align*}
We define the quotient space
\[ \KK\Hom(\underline{\R^{m_\alpha}}, \underline{\R^{n_\beta}}):=\Hom(\underline{\R^{m_\alpha}}, \underline{\R^{n_\beta}})/K_1\cap K_2.\]}\]
defined over $V_\alpha$, such that 
\begin{itemize}
\item[(1.)] $f_1(x)-f_0(x)=\Lambda (s_\alpha(x)), \;\; \forall x\in V_\alpha$;
\item[(2.)] the following diagram defined over $U$  is commutative.
\[\begin{tikzcd}
0 \arrow{r} &\underline{\R^{n_\alpha}} \arrow{d}[swap]{df_1-df_0}\arrow{r}{ds_\alpha} &\underline{\R^{m_\alpha}} \arrow{ld}[swap]{\Lambda}\arrow{d}{\hat{f}_1-\hat{f}_0}\arrow{r} &0 \\
0 \arrow{r} &\underline{\R^{n_\beta}} \arrow{r}{ds_\beta} &\underline{\R^{m_\beta}} \arrow{r} &0
\end{tikzcd}\]
\end{itemize}
Denote this relation by $(f_0,\hat{f}_0)\cong (f_1,\hat{f}_1)$.  
\end{defi}

Recall that in the definition of a smooth manifold, we require to have a collection of local charts $\left\{\phi_i: U_i\ra \R^n\right\}$. The transition maps $\left\{f_{ij}=\phi_j\circ\phi_i^{-1}\right\}$ satisfies the cocycle condition $f_{jk}f_{ij}=f_{ik}$. In the case of Kuranishi manifolds, we use Definition~\ref{def:homotopy} to require the equation $f_{jk}f_{ij}=f_{ik}$ to only hold up to coherent homotopies.

\medskip
\begin{defi}~\label{def:kuranishi}
Let $X$ be a Hausdorff, second countable topological space, and $d\in \Z$.  A Kuranishi atlas on $X$ of virtual dimension $d$ is given by a collection of charts $\Big\{(V_{i},\underline{\R^{m_i}},s_i,\psi_i)\Big\}_{i\in I}$, together with transition morphisms $\Big\{(f_{ij},\hat{f}_{ij})\Big\}$ and homotopy data $\Big\{\Lambda_{ijk}\Big\}$, where
\begin{itemize}
\item[(1.)] $(V_{i},\underline{\R^{m_i}},s_i,\psi_i)$ is a Kuranishi chart with footprint $U_i\subset X$ with $V_i\subset\R^{n_i}$ an open domain, and $n_i-m_i=d$, for all $i\in I$. We require that $\cup_{i\in I} U_i=X$.
\medskip
\item[(2.)] For each $U_{ij}\neq \emptyset$ appearing in the refinement of the open set $U_i\cap U_j$, we have a morphism $(f_{ij}, \hat{f}_{ij}): (V_{\underline{i}j}, \underline{\R^{m_i}}, s_i, \psi_i) \ra (V_{i\underline{j}}, \underline{\R^{m_j}}, s_j, \psi_j)$ of the restricted Kuranishi charts, with $V_{\underline{i}j}\subset V_i$, $V_{i\underline{j}}\subset V_j$ such that $\psi_i^{-1}(V_{\underline{i}j})=U_{ij}$, and $\psi_j^{-1}(V_{i\underline{j}})=U_{ij}$. We also require $(f_{ii}, \hat{f}_{ii})$ to be an open inclusion for all $i\in I$.
\medskip
\item[(3.)] For each $U_{ijk}\neq \emptyset$ appearing in the refinement of $U_{ij}\cap U_{jk}\cap U_{ik}$, there exists open subsets
\[ V_{\underline{i}jk}\subset V_{\underline{i}j}\cap V_{\underline{i}k}, \;\; V_{i\underline{j}k} \subset V_{i\underline{j}}\cap V_{\underline{j}k}, \;\; V_{ij\underline{k}}\subset V_{i\underline{k}}\cap V_{j\underline{k}},\]
with $\psi_i^{-1}(V_{\underline{i}jk})=U_{ijk}$, $\psi_j^{-1}(V_{i\underline{j}k})=U_{ijk}$, $\psi_k^{-1}(V_{ij\underline{k}})=U_{ijk}$, so that the restrictions of $(f_{ij},\hat{f}_{ij})$, $(f_{jk},\hat{f}_{jk})$, $(f_{ik},\hat{f}_{ik})$, and the composition $(f_{jk},\hat{f}_{jk})\circ (f_{ij},\hat{f}_{ij})$ are well-defined. Furthermore, the morphism  $\Lambda_{ijk}:\underline{\R^{m_i}} \ra \underline{\R^{n_k}}$ defined over $V_{\underline{i}jk}$ gives a homotopy $(f_{ik}, \hat{f}_{ik})\cong (f_{jk},\hat{f}_{jk})\circ (f_{ij},\hat{f}_{ij})$
in the sense of Definition~\ref{def:homotopy}. Furthermore, we require that $\Lambda_{iij}=0$ and $\Lambda_{ijj}=0$.
\medskip
\item[(4.)] For each open subset $U_{ijkl}\neq\emptyset$ appearing in the refinement of $U_{ijk}\cap U_{ijl} \cap U_{ikl}\cap U_{jkl}$, we require the equation
\begin{equation}\label{eq:cocycle2}
\Lambda_{ikl}-\Lambda_{jkl}\circ \hat{f}_{ij} -\Lambda_{ijl}+df_{kl}\circ\Lambda_{ijk}=0
\end{equation}
to hold over the topological space $U_{ijkl}$.
\end{itemize}
The topological space $X$ together with a Kuranishi atlas on it is called a Kuranishi manifold.
\end{defi}

\begin{rem}
Equation~\ref{eq:cocycle2} comes from the following diagram:
\[\begin{xy} 
(0,0)*{V_{\underline{i}jkl}}="A"; 
(20,-15)*{V_{i\underline{j}kl}}="B";
(40,0)*{V_{ijk\underline{l}}}="C";
(5,25)*{V_{ij\underline{k}l}}="D";
{\ar@{.>}_{(f_{il},\hat{f}_{il})} "A"; "C"};
{\ar@{->}_{(f_{ij},\hat{f}_{ij})} "A"; "B"};
{\ar@{->}_{(f_{jl},\hat{f}_{jl})} "B"; "C"};
{\ar@{->}^{(f_{ik},\hat{f}_{ik})} "A"; "D"};
{\ar@{->}_{(f_{jk},\hat{f}_{jk})} "B"; "D"};
{\ar@{->}^{(f_{kl},\hat{f}_{kl})} "D"; "C"};
\end{xy}\]
It asserts certain cancellation of the four homotopies on the four facets of the above tetrahedron.
\end{rem}

\medskip
\begin{defi}~\label{def:homo}
Let $\mathfrak{X}$, $\mathfrak{Y}$ be two Kuranishi manifolds, with Kuranishi atlases given by
\begin{align*}\AA &=\big\{ (V_i,\underline{\R^{m_i}}, s_i, \psi_i), (f_{ij},\hat{f}_{ij}), (\Lambda_{ijk})\big\}_{i,j,k\in I}, \\
\BB &= \big\{ (V_p,\underline{\R^{m_p}}, s_p, \psi_p), (f_{pq},\hat{f}_{pq}), (\Lambda_{pqr}) \big\}_{p,q,r\in P}\end{align*}
A {\sl strict Kuranishi morphism} $h: \mathfrak{X}\ra\mathfrak{Y}$ is given by the following data:
\begin{itemize}
\item[(1.)] A map of indices $\tau: I \ra P$.
\item[(2.)] For each $i\in I$, a morphism
\[ (h_{i}, \hat{h}_{i}): (V_{i}, \underline{\R^{m_{i}}}, s_{i}, \psi_{i}) \ra (V_{\tau(i)}, \underline{\R^{m_{\tau(i)}}}, s_{\tau(i)}, \psi_{\tau(i)})\]
between Kuranishi charts.
\item[(3.)] For each pair of indices $(i,j)\in I\times I$ such that $U_{ij}\neq \emptyset$, a quotient bundle map $\Delta_{ij}\in\KK\Hom(\underline{\R^{m_i}}, \underline{\R^{n_{\tau(j)}}})$, giving a homotopy
\[ \Delta_{ij}: (h_{j},\hat{h}_{j})\circ f_{ij} \cong f_{\tau(i)\tau(j)}\circ (h_{i},\hat{h}_{i}). \]
\item[(4.)] For a triple of indices $(i,j,k)\in I\times I\times I$ such that $U_{ijk}\neq \emptyset$, we require the equation
\begin{equation}~\label{eq:2-hom} \Delta_{ik}-dh_{k}\circ\Lambda_{ijk}+\Lambda_{\tau(i)\tau(j)\tau(k)}\circ \hat{h}_{i}-df_{\tau(j)\tau(k)}\circ\Delta_{ij}-\Delta_{jk}\circ \hat{f}_{ij}=0\end{equation}
to hold over $U_{ijk}$.
\end{itemize}
\end{defi}
\begin{rem}
Again, the appearance of Equation~\ref{eq:2-hom} asserts certain cancellations of homotopies. This is necessary to define the tangent complexes associated with a morphism, see Definition~\ref{def:tangent}.
\end{rem}

In a similar fashion, one may also define $2$-morphisms between strict Kuranishi morphisms. 

\medskip
\begin{defi}~\label{def:equihom}
Let $$h, g: (X,\AA) \ra (Y,\BB)$$ be two strict Kuranishi morphisms between Kuranishi manifolds, both subjected to the atlases $\AA$ and $\BB$ as in Definition~\ref{def:homo}. We say that $h$ is homotopic to $g$ if the following conditions hold.
\begin{itemize}
\item[(a.)] The underlying continuous maps $\underline{h}=\underline{g}$.
\item[(b.)] For each $i\in I$, the following diagram is commutative, up to homotopy.
\[\begin{xy} 
(0,0)*{h_i^{-1}(V_{\underline{\tau_h(i)}\tau_g(i)})\cap g_i^{-1}(V_{\tau_h(i)\underline{\tau_g(i)}})}="A"; (80,0)*{V_{\tau_h(i)\underline{\tau_g(i)}}}="C"; (40,20)*{V_{\underline{\tau_h(i)}\tau_g(i)}}="B";
{\ar@{->}^{(g_i,\hat{g}_i)} "A"; "C"};
{\ar@{->}|-{(h_i,\hat{h}_i)} "A"; "B"};
{\ar@{->}|-{(f_{\tau_h(i)\tau_g(i)},\hat{f}_{\tau_h(i)\tau_g(i)})} "B"; "C"};
\end{xy}\]
In other words, there exists a homotopy $\Upsilon_i: \underline{\R^{m_i}} \ra \underline{\R^{m_{\tau_g(i)}}}$ between morphisms $(g_i, \hat{g}_i)$ and $(f_{\tau_h(i)\tau_g(i)}, \hat{f}_{\tau_h(i)\tau_g(i)}) \circ (h_i,\hat{h}_i)$.
\item[(c.)] For a pair of indices $(i,j)\in I\times I$, we require that
\[\Delta^g_{ij}-df_{\tau_g(i)\tau_g(j)} \circ \Upsilon_i + \Upsilon_j\circ \hat{f}_{ij} - df_{\tau_h(j)\tau_g(j)} \circ \Delta^h_{ij} - [\Lambda_{\tau_h(i)\tau_g(i)\tau_g(j)}+\Lambda_{\tau_h(i)\tau_h(j)\tau_g(j)}]\circ \hat{h}_i=0,
\]
as morphisms over the intersection $U_{ij}$.
\end{itemize}
\end{defi}

In a forthcoming paper~\cite{Tu-1}, we prove the following 

\medskip
\begin{thm}~\label{conj-main}
There exists a strict $2$-category $\mathfrak{Kur}$ whose objects are Kuranishi manifolds, with the following properties:
\begin{itemize}
\item[(A.)] The category of smooth manifolds $\mathfrak{Man}$ imbeds fully faithfully into $\mathfrak{Kur}$.
\item[(B.)] There exists a $2$-categorical fiber product $\mathfrak{X}\times_M\mathfrak{Y}$ in $\mathfrak{Kur}$ if the base $M$ is a smooth manifold (which makes sense due to part $(A.)$).
\[\begin{CD}
\mathfrak{X}\times_M\mathfrak{Y} @>>> \mathfrak{Y}\\
@VVV                  @VVV\\
\mathfrak{X}@>>> M.
\end{CD}\]
\end{itemize}
\end{thm}

\medskip
\noindent {\sl Sketch of proof.} Putting the definitions of Kuranishi manifolds, strict morphisms, and $2$-morphisms together, one can prove that we obtain a strict $2$-category ${\sf pre-}\mathfrak{Kur}$. This construction is rather straightforward, and such a structure essentially follows from the fact that the category of $2$-term complexes of vector bundles on a topological space forms a strict $2$-category.  

However, the category ${\sf pre-}\mathfrak{Kur}$ is not what we want, as the category of manifolds $\mathfrak{Man}$ does not imbeds into ${\sf pre-}\mathfrak{Kur}$. The problem is that when writing down a smooth map between two manifolds in terms atlases, one should be able to allow refinements on the domain manifolds' atlas. The solution is to single out a special class $\mathfrak{R}$ of $1$-morphisms in ${\sf pre-}\mathfrak{Kur}$ consisting of refinements of Kuranshi atlases. The desired category should then be given by
\[ \mathfrak{Kur}:= \mathfrak{R}^{-1}\big({\sf pre-}\mathfrak{Kur}\big),\]
the $2$-categorical localization of ${\sf pre-}\mathfrak{Kur}$ with respect to $\mathfrak{R}$. There is a complication involved here when forming a $2$-categorical localization. If we were only localizing an ordinary category, to have a reasonable localization, one only needs the class $\mathfrak{R}$ be to a multiplicative system. While it is indeed the case that $\mathfrak{R}$ is multiplicative, it does not guarantee the existence of a $2$-categorical localization. However, one can show that $\mathfrak{R}$ has the important property that there exists a {\sl canonical} base change for morphisms in $\mathfrak{R}$, illustrated in the following (strictly commutative) diagram.
\[\begin{CD}
\mathfrak{X}'  @>\exists f'>> \mathfrak{Y}'\\
@V \exists s\in\mathfrak{R} VV @V r\in \mathfrak{R} VV\\
\mathfrak{X} @>f>> \mathfrak{Y}.
\end{CD}\]
With this extra property, it is not difficult to prove that the ordinary localized category $\mathfrak{Kur}=\mathfrak{R}^{-1}\big({\sf pre-}\mathfrak{Kur}\big)$ with morphisms given by roof diagrams, admits a natural $2$-category structure. This proves part $(A.)$ of the theorem. Part $(B.)$ is lengthy, but rather straightforward.\ed

\begin{rem}
Part $(A.)$ of the above theorem should be attributed to Joyce~\cite{Joyce}, and a more general version of part $(B.)$ was also announced in {\sl loc. cit.}. Strictly speaking, the definition of Kuranishi manifolds by Joyce is different from Definition~\ref{def:kuranishi}. However, we expect the homotopy category of the $2$-category $\mathfrak{Kur}$ to be equivalent to Joyce's version. The main point we want to make here is that $[0,1]$-type homotopy $L_\infty$ spaces are automatically Kuranishi manifolds by its very formulation. Indeed, from the $L_\infty$ point of view, Joyce's definition of coordinate changes is quite natural.
\end{rem}

\section{Virtual fundamental cycles}

In this section, we consider the fundamental problem of constructing fundamental cycles on $[0,1]$-type homotopy $L_\infty$ spaces. Since these spaces are automatically Kuranishi manifolds, it suffices to perform the construction for Kuranishi manifolds. More precisely, in a forthcoming paper~\cite{Tu-2}, we prove the following

\medskip
\begin{thm}~\label{main-conj2}
Let $\mathfrak{X}$ be an oriented Kuranishi manifold of virtual dimension $d$. Assume that its underlying topological space $X$ is a finite CW-complex. Then there exist natural constructions of a homology class $[\mathfrak{X}]_{\virt}\in H_d(X)$ and an oriented cobordism class $[\mathfrak{X}]_{\virt}^{\cob}\in \Omega^{\sf or}_d$. Furthermore, the virtual cobordism class is deformation invariant in the following sense. If $\mathfrak{X}_t$ is a family of Kuranishi manifolds parametrized by $t\in [0,1]$, and given a family of fiber diagrams (with $N$ a manifold)
\[\begin{CD}
\mathfrak{X}_t \times_N \mathfrak{Y} @>>> \mathfrak{Y}\\
@VVV @VVV\\
\mathfrak{X}_t @>>> N,
\end{CD}\]
we have the class $[\mathfrak{X}_t \times_N \mathfrak{Y}]_\virt^{\cob}\in \Omega^{\sf or}_d$ is independent of $t$.
\end{thm}

\medskip
\begin{rem}~\label{rem-def}
Two special cases are:
\begin{itemize}
\item[--] Let $\mathfrak{Y}$ and $N$ both be a point. We get that the class $[\mathfrak{X}_t]_{\virt}^{\cob}$ is independent of $t$.
\item[--] Let $\mathfrak{X}_0$ and $\mathfrak{Y}$ both be submanifolds of the manifold $N$, of complimentary dimensions. We obtain that the number $[\mathfrak{X}_0\times_N \mathfrak{Y}]_\virt^{\cob}\in \Omega^{\sf or}_0\cong \Z$ is given by the intersection number between $\mathfrak{X}_0$ and $\mathfrak{Y}$ in $N$.
\end{itemize}
\end{rem}

The main scheme of the proof of the above theorem follows Spivak~\cite{Spivak}. It uses the construction of a tangent complex associated with a Kuranishi morphism which we describe now. Let $h: \mathfrak{X} \ra \mathfrak{Y}$ be a strict morphism of Kuranishi manifolds. We continue to use notations as in Definition~\ref{def:homo}. For each $i\in I$, the morphism $(h_i,\hat{h}_i)$ induces a morphism between $2$-term complexes:
\[ (dh_i,\hat{h}_i): \TT_i \ra \TT_{\tau(i)},\;\; \TT_i:=\underline{\R^{n_i}} \ra \underline{\R^{m_i}}, \;\;\TT_{\tau(i)}:=\underline{\R^{n_{\tau(i)}}} \ra \underline{\R^{m_{\tau(i)}}}.\]
Consider its mapping cone
\[ \Cone(dh_i,\hat{h}_i):= \TT_i \oplus \TT_{\tau(i)} [-1], \;\; \]
or explicitly the complex
\[\begin{CD}
0 @>>> \underline{\R^{n_i}} @>\begin{bmatrix} ds_i \\ dh_i \end{bmatrix}>> \begin{matrix}\underline{\R^{m_i}}\\ \oplus\\ \underline{\R^{n_{\tau(i)}}}\end{matrix} @>\big[\hat{h}_i \; -ds_{\tau(i)}\big]>> \underline{\R^{m_{\tau(i)}}} @>>>0.
\end{CD}\]
For a pair of indices $(i,j)$ such that $U_{ij}\neq \emptyset$, the existence of a homotopy
$\Delta_{ij}$ implies there exists a morphism of between the mapping cones over the intersection $U_{ij}$, explicitly described in the following diagram.

\medskip
\begin{equation}\label{eq:trans}\begin{CD}
0 @>>> \underline{\R^{n_i}} @>\begin{bmatrix} ds_i \\ dh_i \end{bmatrix}>> \begin{matrix}\underline{\R^{m_i}}\\ \oplus\\ \underline{\R^{n_{\tau(i)}}}\end{matrix} @>\big[\hat{h}_i \; -ds_{\tau(i)}\big]>> \underline{\R^{m_{\tau(i)}}} @>>>0.\\
@. @V df_{ij} VV @VV\begin{bmatrix} \hat{f}_{ij} \;\;\;\;\;\;\;\;\;\;\;\;\;\; 0 \\ -\Delta_{ij} \;\; df_{\tau(i)\tau(j)}\end{bmatrix} V @VV \hat{f}_{\tau(i)\tau(j)} V @.\\
0 @>>> \underline{\R^{n_j}} @>\begin{bmatrix} ds_j \\ dh_j\end{bmatrix}>> \begin{matrix}\underline{\R^{m_j}}\\ \oplus\\ \underline{\R^{n_{\tau(j)}}}\end{matrix} @>\big[\hat{h}_j \; -ds_{\tau(j)}\big]>> \underline{\R^{m_{\tau(j)}}} @>>>0.
\end{CD}\end{equation}
\begin{lem}
The collection of morphisms $$\Big(df_{ij}, \begin{bmatrix} \hat{f}_{ij} \;\;\;\;\;\;\;\;\;\;\;\;\;\; 0 \\ -\Delta_{ij} \;\; df_{\tau(i)\tau(j)}\end{bmatrix}, \hat{f}_{\tau(i)\tau(j)}\Big), \;\; (i,j)\in I\times I$$ between the mapping cones satisfies weak cocycle condition on triple intersections.
\end{lem}

\medskip
\begin{defi}~\label{def:tangent}
Let $h: \mathfrak{X} \ra \mathfrak{Y}$ be a morphism between two Kuranishi manifolds. For each $i\in I$, set
\[ T_{h,i}^*:=\underline{H}^*(\Cone(dh_i,\hat{h}_i)), \;\; *=0,1,2.\]
For each $*=0,1,2$, the previous lemma implies that the collection $\left\{T_{h,i}^*\right\}_{i\in \II}$ satisfies the strict cocycle condition. Thus we obtain three globally defined sheaf on $X$, denoted by $T_h^*\; (*=0,1,2)$. We refer to them as the tangent sheaves of the morphism $h$. 
\end{defi}

In manifold theory, the  Whitney's imbedding theorem states that a compact manifold can always be realized as a submanifold of an Euclidean space $\R^N$ for some large $N$. The construction of such an imbedding in the manifold case can be used to prove the following result.

\begin{prop}
Let $\mathfrak{X}$ be a compact Kuranishi manifold of virtual dimension ${\sf vdim}\mathfrak{X}$. Then there exists a morphism 
\[ h: \mathfrak{X}\ra \R^N\]
for some large enough $N$, such that
\begin{itemize}
\item[(a.)] the tangent sheaves $T^0_h$, $T^2_h$ are zero.
\item[(b.)] the tangent sheaf $T_h^1$ is a vector bundle on $X$ of rank $N -{\sf vdim} \mathfrak{X}$.
\item[(c.)] the underlying map of topological spaces $\underline{h}: X \ra \R^N$ is a homeomorphism onto its image.
\end{itemize}
\end{prop}

We proceed to construct the virtual fundamental cycle of $\mathfrak{X}$. For this, we make the following

\medskip
\noindent {\bf Assumption $(\star)$.} The Kuranishi manifold $\mathfrak{X}$ is oriented, and its underlying topological space $X$ is a finite (hence compact) CW-complex.

\medskip
Now, let $\mathfrak{X}$ be a Kuranishi manifold satisfying Assumption $(\star)$. In particular, it is compact. The proposition above implies a Kuranishi morphism $h: \mathfrak{X} \ra \R^N$ with conditions $(a.),(b.),(c.)$. Since $X$ is a finite CW-complex, its image $h(X)\subset \R^N$ is an Euclidean Neighborhood Retract (ENR). We choose such an open neighborhood $W\supset h(X)$, together with a retract $$p: W\ra X$$ (so $p\circ h=\id$). For a given $\epsilon>0$, we set 
\[ W_\epsilon:= \left\{ w\in W \mid ||w-p(w)||<\epsilon.\right\}\]
The restriction of $p$ to $W_\epsilon$ is still a retraction which we continue to denote by $p$. Using a partition of unity argument, one can prove the following

\medskip
\begin{prop}
Let $\mathfrak{X}$ be a Kuranishi manifold satisfying Assumption $(\star)$. Then there exists a small enough $\epsilon>0$, and a section $s$ of the bundle $p^*T_h^1$ over $W_\epsilon$, such that $s^{-1}(0)=h(X)$. Furthermore, there exists a small perturbation $s'$ of the section $s$ such that  $s'^{-1}(0)$ is a compact oriented manifold of dimension ${\sf vdim} \mathfrak{X}$. 
\end{prop}

\medskip
\begin{defi}~\label{def:vfc}
Let $\mathfrak{X}$ be an oriented Kuranishi manifold of virtual dimension $d$. Assume that its underlying topological space $X$ a finite CW-complex. Then, its virtual fundamental class is defined to be
\[ [\mathfrak{X}]_{\virt}:= p_*[s'^{-1}(0)]\in H_{\sf vdim}(X).\]
Its virtual cobordism class is defined to be $[\mathfrak{X}]_\virt^{\cob}:=[s'^{-1}(0)]\in \Omega^{\sf or}_d$. 
\end{defi}

Implicit in the above definition is that the homology class $[\mathfrak{X}]_{\virt}$ (or $[\mathfrak{X}]_\virt^{\cob}$) is independent of all choices involved in this construction. Indeed, the choices include the morphism $h$, the retraction $p$, the constant $\epsilon$, the section $s$ (where a partition of unity is chosen on $W_\epsilon$), and the perturbation $s'$. Proving such independence is non-trivial, but is of technical nature. The moral is that all the choices made can vary in a continuous family, hence inducing homotopic cycles.

Back to Casson invariants, let $M$ be an oriented $3$-dimensional integral homology sphere. By Theorem~\ref{thm-lie},  the moduli space $R^\irr(M)$ admits a natural oriented $[0,1]$-type homotopy $L_\infty$ enhancement $\mathfrak{R}^\irr(M)$ of virtual dimension $0$. This implies in particular that $\mathfrak{R}^\irr(M)$ is also an oriented Kuranishi manifold  of virtual dimension $0$ whose underlying topological space is $R^\irr(M)$. We also note that the compact topological space $R^\irr(M)$ is a real analytic space by Proposition~\ref{prop}, which implies that it is a finite CW complex~\cite[Section 4]{VM}. Thus applying Theorem~\ref{conj-main}, we obtain a homology class
\[ [\mathfrak{R}^\irr(M)]_\virt\in H_0(R^\irr(M)).\]

\medskip
\begin{thm}~\label{conj-casson}
The Casson invariant $\lambda(M)$ is given by
\begin{equation}~\label{eq-main} \lambda(M)=\int_{[\mathfrak{R}^\irr(M)]_{\virt}} 1/2.\end{equation}
\end{thm}

\medskip
\noindent {\sl Sketch of proof.} Let $M=H_1\cup_\Sigma H_2$ be the Heegaard splitting used in the definition of Casson invariant of $M$. It induces a covering of $M$ by two open subsets $M_1$ and $M_2$, with $H_1$ a deformation retract of $M_1$, $H_2$ a deformation retract of $M_2$, and $\Sigma$ a deformation retract of $M_1\cap M_2$. The \u Cech complex associated to this covering yields a fiber diagram of differential graded Lie algebras
\[ \begin{CD}
\Omega^*(M)\otimes \su(2) @>>> \Omega^*(M_1)\otimes\su(2)\\
@VVV               @VVV\\
\Omega^*(M_2)\otimes \su(2) @>>> \Omega^*(M_1\cap M_2)\otimes \su(2).
\end{CD}\]
Theorem~\ref{thm-lie} applied to the above diagram, together with the invariance of de Rham cohomology under deformation retracts, gives a commutative diagram of Kuranishi manifolds
\begin{equation}~\label{diagram:kuranishi}\begin{CD}
\mathfrak{R}^\irr(M)   @>>> \mathfrak{R}^\irr(H_1)\\
@VVV   @VVV\\
\mathfrak{R}^\irr(H_2) @>>> \mathfrak{R}^\irr(\Sigma),
\end{CD}\end{equation}
enhancing Diagram~\ref{diagram:classical}. Now, since both $H_1$ and $H_2$ retracts to wedge products of circles, there are no cohomology in degrees greater or equal to $2$, which implies that the Kuranishi manifolds $\mathfrak{R}^\irr(H_1)$ and $\mathfrak{R}^\irr(H_2)$ agree with the ordinary manifolds $R^\irr(H_1)$ and $R^\irr(H_2)$. For the Kuranishi $\mathfrak{R}^\irr(\Sigma)$, by Poincar\'{e} duality the triviality of degree $0$ cohomology implies the triviality of degree $2$ cohomology, which proves that $\mathfrak{R}^\irr(\Sigma)$ is simply the underlying manifold $R^\irr(\Sigma)$.

Furthermore, the fact that the \u Cech complex is a resolution of the de Rham complex of $M$ makes Diagram~\ref{diagram:kuranishi} above a fiber diagram in the $2$-category $\mathfrak{Kur}$. To this end, we take a small perturbation $\mathfrak{R}^\irr(H_2)_t \; (t\in [0,1])$ of the submanifold $\mathfrak{R}^\irr(H_2)$. At $t=0$ where we do not perturb, the fiber product is the Kuranishi manifold $\mathfrak{R}^\irr(M)$, while at $t=1$ it gives the transversal intersection $\mathfrak{R}^\irr(H_2)_1 \cap \mathfrak{R}^\irr(H_1)$. This proves the theorem by deformation invariance of the virtual fundamental class (see Theorem~\ref{main-conj2} and Remark~\ref{rem-def}).\ed

\end{document}